\documentclass[11pt]{amsart}

\usepackage{amssymb,amscd,amsthm,amsxtra}
\usepackage{latexsym}

\newtheorem{thm}{Theorem}[section]

\newtheorem{lem}[thm]{Lemma}
\newtheorem{prop}[thm]{Proposition}
\theoremstyle{definition}
\newcommand{\comment}[1]{}

\newtheorem{defn}[thm]{Definition}
\theoremstyle{remark}
\newtheorem{rem}[thm]{Remark}
\numberwithin{equation}{section}
\newcommand{\R}{\mathbb R}
\newcommand{\C}{\mathbb C}

\begin{document}

\title[Global Well-Posedness for the $L^2$-critical nonlinear Schr\"odinger equation]
{\bf Global Well-Posedness for the $L^2$-critical nonlinear
Schr\"odinger equation in higher dimensions}

\author{Daniela De Silva}
\address{Department of Mathematics, Johns Hopkins University, Baltimore, MD 21218}
\email{\tt  desilva@math.jhu.edu}

\author{Nata\v{s}a Pavlovi\'{c}}
\address{Department of Mathematics, Princeton University, Princeton, NJ 08544-1000}
\email{\tt natasa@math.princeton.edu}

\author{Gigliola Staffilani}
\address{Department of Mathematics, Massachusetts Institute of Technology,
Cambridge, MA 02139-4307}
\email{\tt gigliola@math.mit.edu}

\author{Nikolaos Tzirakis}
\address{Department of Mathematics, University of Toronto, Toronto, Ontario, Canada M5S 2E4}
\email{tzirakis@math.toronto.edu}
\thanks{N.P. was supported by N.S.F. Grant DMS 0304594 and
G.S. was supported by N.S.F. Grant DMS 0602678.}
\date{July 24, 2006}


\begin{abstract}
The initial value problem for the  $L^{2}$ critical semilinear
Schr\"odinger equation in $\R^n, n \geq 3$  is considered. We show
that the problem is globally well posed in $H^{s}({\Bbb R^{n} })$
when $1>s>\frac{\sqrt{7}-1}{3}$ for $n=3$, and when $1>s>
\frac{-(n-2)+\sqrt{(n-2)^2+8(n-2)}}{4}$ for $n \geq 4$. We use the
``$I$-method'' combined with a local in time Morawetz estimate.
\end{abstract}
\maketitle

\section{Introduction}

In this paper we study the $L^{2}$ critical defocusing Cauchy problem

\begin{align}
&iu_{t}+ \Delta u -|u|^{\frac{4}{n}}u=0, \; \; x \in \R^{n}, \; t \geq 0, \label{ivp1}\\
&u(x,0)=u_{0}(x)\in H^{s}(\R^n), \; \; n\geq 3,\label{bc1}
\end{align}
where $u(t,x)$ is a complex-valued function in space-time $\R^+
\times \R^n$. Here $H^s(\R^n)$ denotes the usual inhomogeneous
Sobolev space.

The local well-posedness definition that we use here reads as
follow: for any choice of initial data $u_0 \in H^s$, there exists
a positive time $T = T(\|u_0\|_{H^{s}})$ depending only on the
norm of the initial data, such that a solution to the initial
value problem exists on the time interval $[0,T]$, it is unique in
a certain Banach space of functions $X\subset C([0,T],H^s_{x})$,
and the solution map from $H^s_{x}$ to $C([0,T],H^s_{x})$ depends
continuously on the initial data on the time interval $[0,T]$. If
$T = \infty$ we say that a Cauchy problem is globally well-posed.

Although here we only study the case $n \geq 3$, we recall  some
known results for general dimensions that highlight, in a certain
way, the differences that arise when the nonlinearity is not
smooth. It is known, for example, that the local and global theory
for \eqref{ivp1}-\eqref{bc1}, if one considers general smooth
data, depend on the smoothness of the nonlinearity in a crucial
way. In the case when the nonlinearity is smooth enough (for
example, for $n=1$ and $n=2$, when the nonlinearity is algebraic
and thus $C^{\infty}$) regularity properties of solutions to the
above initial value problem are very well understood. However in
the general case, certain restrictions on $s$ are needed in order
to answer the questions of local/global well-posedness, regularity
and others, as clearly as in the algebraic case. For more
information the reader should consult \cite{cw1}.

For our purposes we restrict ourselves to initial data in $H^{s}$
with $0<s<1$, and make some comments for the limiting cases where
$s=0, 1$. It is known, \cite{cw}, \cite{cw1}, \cite{tk1}, that the
initial value problem \eqref{ivp1}-\eqref{bc1} is locally well-posed
in $H^s$ when $s>0,$ and local in time solutions enjoy mass
conservation
\begin{equation}\label{mass}
\|u(\cdot,t)\|_{L^2(\R^n)} =
\|u_0(\cdot)\|_{L^2(\R^n)}.
\end{equation}
Moreover, $H^{1}$ solutions enjoy conservation of the energy
\begin{equation}
E(u)(t)=\frac{1}{2}\int |\nabla u(t)|^{2}dx+\frac{n}{2n+4}\int
|u(t)|^{\frac{2n+4}{n}}dx=E(u)(0),
\end{equation} which together with \eqref{mass} and the local theory immediately
yields global in time well-posedness for \eqref{ivp1}-\eqref{bc1}
with initial data in $H^1.$ Actually T. Kato proved local and
global well posedness at the energy level for nonlinearities that
are more general than ours, \cite{tk1}. Also in \cite{cw}, it is
proved that for data in $L^{2}$ the IVP \eqref{ivp1}-\eqref{bc1} is
well-posed in an interval $[0,T]$, but in this case $T=T(u_{0})$,
making the conservation of the mass not of immediate use in order to
obtain global well-posedness.

The main purpose of this paper is to partially extend the
techniques that have been developed so far to obtain local or
global well-posedness for the $L^{2}$ critical problem to the non
algebraic case.

The question of global well-posedness of
\eqref{ivp1}-\eqref{bc1} in the case $0 \leq s < 1$  is at this
point only partially answered. Previous work establishes that in one
 dimension  the problem is globally well-posed for $s>4/9$
(see \cite{ckstt3}, \cite{tz}) and  in two dimensions for $s \geq 1/2$
(see \cite{ckstt5}  \cite{FG}).
In both these cases no scattering has been proved.
%In three dimensions the problem is globally well-posed
%for $s>4/5$ and a scattering result is also available\footnote{This is
%because the proof is based on an  \emph{a priori} global space-time bound for the
%solution obtained through an interaction Morawetz estimate.}, \cite{ckstt6}.
In higher dimensions a first result on global well-posedness below
the energy norm was presented in \cite{ckstt7}. There the authors
obtained $s>1-\epsilon$, with $\epsilon>0$, but not explicitly
quantifiable. Recently T. Tao and M. Visan announced that
\eqref{ivp1}-\eqref{bc1} is globally well-posed in $L^{2}$
provided $n=4$ and the data are radially symmetric.  Still in
higher dimensions a partial result is included in \cite{vz2}
without being explicitly stated. In fact there the authors are
considering the $L^{2}$ critical focusing problem. However, a
byproduct of their analysis using the ``$I$-method'' gives global
well-posedness for some $s<1$ in dimensions $n \geq 3$. In this
paper we extend this result. The precise statement of what we
prove is contained in the following theorem.

\begin{thm}\label{main}
The initial value problem \eqref{ivp1}-\eqref{bc1} is globally
well-posed in $H^{s}(\R^n)$, for any $1>s>\frac{\sqrt{7}-1}{3}$ when
$n=3$, and for any $1>s
> \frac{-(n-2)+\sqrt{(n-2)^2+8(n-2)}}{4}$ when $n \geq 4.$

\end{thm}

We notice that our best result, which is obtained for $n=3$, gives
$s>0.55$ and as $n \rightarrow \infty$, $s \rightarrow 1$.

The proof of Theorem \ref{main} was inspired by a recent paper of
Fang and Grillakis \cite{FG} in which the case $n=2$ is
considered. It relies on two main ingredients. The first one is
the so-called $I$-method introduced  by J. Colliander, M. Keel, G.
Staffilani, H. Takaoka and T. Tao (see, for example \cite{ckstt1,
ckstt5, ckstt4, ckstt2}) which is based on the almost conservation
of a certain modified energy functional. The idea is to replace
the conserved quantity $E(u)$ which is no longer available for
$s<1$, with an ``almost conserved'' variant $E(Iu)$ where $I$ is a
smoothing operator of order $1-s$ which behaves like the identity
for low frequencies and like a fractional integral operator for
high frequencies. Thus the operator $I$ maps $H_{x}^{s}$ to
$H_{x}^{1}$. Notice that $Iu$ is not a solution to \eqref{ivp1}
and hence we expect an energy increment. This increment is in fact
quantifying $E(Iu)$ as an ``almost conserved'' energy. The key is
to prove that on intervals of fixed length, where local
well-posedness is satisfied, the increment of the modified energy
$E(Iu)$ decays with respect to a large parameter $N$. (For the
precise definition of $I$ and $N$ we refer the reader to Section
$2$.) This requires delicate estimates on the commutator between
$I$ and the nonlinearity. In dimensions 1 and 2, where the
nonlinearity is algebraic, one can write the commutator explicitly
using the Fourier transform, and control it by multi-linear
analysis and bilinear estimates. The analysis above can be carried
out in the $X^{s,b}$ spaces setting, where one can use the
smoothing bilinear Strichartz  estimate of Bourgain (see e.g.
\cite{bo3}) along with Strichartz estimates to demonstrate the
existence of global rough solutions (see \cite{ckstt1, ckstt3,
ckstt5} and \cite{tz}). In our case,  where $n \geq 3$, we cannot
use any sort of smoothing estimates and we rely purely on
Strichartz estimates alone. One of the advantages of the
``$I$-method'' is that one can use commutator estimates involving
the operator $I$ as a substitute for smoothing estimates even when
the nonlinearity has no smoothing properties.

The second ingredient in our proof is an \textit{a priori}
interaction Morawetz-type estimate for the solution $u$ to
\eqref{ivp1} (see \cite{ckstt6}, \cite{FG}, \cite{tvz}). With the
help of this estimate, some harmonic analysis and interpolation we
obtain for any compact interval $J=[a,b]$ a priori control of the
$L^{\frac{4(n-1)}{n}}_tL^{\frac{2(n-1)}{n-2}}_x(J\times \R^n )$
mixed type Lebesgue norm. More precisely we have,
\begin{equation}\label{mor}
\|u\|_{L^{\frac{4(n-1)}{n}}_tL^{\frac{2(n-1)}{n-2}}_x(J\times \R^n
)} \lesssim (b-a)^{\frac{n-2}{4(n-1)}}
\|u_0\|_{L^2_x}^{\frac{1}{2}}\|u\|_{L^\infty_t\dot{H}^{\frac{1}{2}}_{x}(J\times
\R^n )}^{\frac{n-2}{n-1}}.
\end{equation}
Notice that the above norm is Strichartz admissible (again consult
Section $2$ for a definition) a fact that will be very important
in our argument.

Combining the two ingredients above, the idea of the proof is as
follows. Fix a large value of time $T_{0}$. We observe that if $u$
is a solution to \eqref{ivp1} in the time interval $[0,T_{0}],$
then
$u^\lambda(x)=\frac{1}{\lambda^{\frac{n}{2}}}u(\frac{x}{\lambda},\frac{t}{\lambda^{2}})$
is a solution to the same equation in $[0, \lambda^{2}T_{0}]$. We
choose the parameter $\lambda>0$ so that
$E(Iu_{0}^{\lambda})=O(1)$. Using Strichartz estimates we show
(see Proposition 3.3) that if $J=[a,b]$ and
$\|u^{\lambda}\|_{L^{\frac{4(n-1)}{n}}_tL^{\frac{2(n-1)}{n-2}}_x(J\times
\R^n )}^{\frac{4(n-1)}{n}}<\mu$, where $\mu$ is a small universal
constant, then
$$\label{ZI} Z_I(J) :=
\sup_{(q,r) \ \ admissible} \|\langle \nabla \rangle
Iu^{\lambda}\|_{L^q_tL^r_x(J\times \R^n)} \lesssim \|I
u^\lambda(a)\|_{H^1}.$$ Moreover in  this same time interval where
the problem is then well-posed, we can prove the ``almost
conservation law'' (see Proposition 4.1)
\begin{equation}\label{acl}
    |E(Iu^{\lambda})(b) - E(Iu^{\lambda})(a)| \lesssim N^{-1+s-\nu}
Z_I(J)^{2+\frac{8}{n}} \lesssim N^{-1+s-\nu}\|I
u^\lambda(a)\|^{2+\frac{8}{n}}_{H^1},\end{equation} for some
suitable $\nu>0$.

Of course for the arbitrarily large
interval $[0,\lambda^{2}T_{0}] $ we do not have that
$$\|u^{\lambda}\|_{L^{\frac{4(n-1)}{n}}_tL^{\frac{2(n-1)}{n-2}}_x([0,\lambda^{2}T_{0}]\times
\R^n )}^{\frac{4(n-1)}{n}}<\mu.$$
This is where we use (\ref{mor}). We first
control the growth of (\ref{mor}) in $[0,\lambda^{2}T_{0}]$.
A little analysis shows that
\begin{align}
& \|u^\lambda\|_{L^{\frac{4(n-1)}{n}}_tL^{\frac{2(n-1)}{n-2}}_x([0,\lambda^{2}T_{0}]\times
\R^n )} \nonumber \\
& \lesssim (\lambda^{2}T_{0})^{\frac{n-2}{4(n-1)}}
\|u_0\|_{L^2_x}^{\frac{1}{2}}\sup_{[0,\lambda^{2}T_{0}]}\left(\|
u_0\|^{\frac{1}{2}}_{L^2}\| I u^\lambda(t)\|^{\frac{1}{2}}_{H^1} +
\|I u^\lambda(t)\|_{H^1}\right)^{\frac{n-2}{n-1}}. \label{intMorfinal3}
\end{align}
Now suppose  we knew that for any $t\in [0,\lambda^{2}T_{0}]$,
\begin{equation}\label{persist}
    \|Iu^\lambda(t)\|_{H^1}=O(1).
    \end{equation}
Then we partition the  arbitrarily large interval
$[0,\lambda^{2}T_{0}]$ into $L$ intervals where the local theory
uniformly applies.  From \eqref{intMorfinal3} we have
\begin{equation} \label{intL}
L \sim \frac{C(\lambda^{2}T_{0})^{\frac{n-2}{n}}}{\mu},
\end{equation}
with $C$ a large constant that will depend only on the
$\|u_{0}\|_{L^{2}}$ norm. $L$  is of course
finite and defines  the number of the intervals in the partition
that will make the Strichartz
$L^{\frac{4(n-1)}{n}}_tL^{\frac{2(n-1)}{n-2}}_x$ norm less than
$\mu$.

In order to obtain \eqref{persist} we observe that
\[\|Iu^\lambda(t)\|_{H^1}\leq
E^{1/2}(Iu^\lambda(t))+\|u^\lambda(t)\|_{L^{2}}.\]
Hence by the fact that the $L^{2}$ norm is scaling invariant and
conserved, we only have to show that $E(Iu^\lambda(t))=O(1)$ for
all $t\in [0,\lambda^{2}T_{0}]$. By the ``almost conservation''
law \eqref{acl} we then require that  $L \sim N^{1-s+\nu}$, for
suitable $\nu$. Since this restriction needs to be compatible with
\eqref{intL}, we obtain the conditions $1>s>\frac{\sqrt{7}-1}{3}$
for $n=3$ and $1>s>\frac{-(n-2)+\sqrt{(n-2)^2+8(n-2)}}{4}$ for any
$n \geq 4$. For a more detailed proof the reader should check
Section 5.
\\
\\
The rest of the paper is organized as follows. In Section 2 we
introduce some notation and state important propositions that we
will use throughout the paper. There we also present as in
\cite{ckstt6},\cite{tvz} the estimate (\ref{mor}). In Section 3 we
prove the local well-posedness theory for $Iu,$ and the main
estimates that we use to prove the decay of the increment of the
modified energy. The decay itself is obtained in Section 4.
Finally in Section 5 we give the details of the proof of global
well-posedness stated in Theorem \ref{main}.

\section{Preliminaries}

In what follows we use $A \lesssim B$ to denote an estimate of the form $A\leq CB$ for some
absolute constant $C$.
If $A \lesssim B$ and $B \lesssim A$ we say that $A \sim B$. We write $A \ll B$ to denote
an estimate of the form $A \leq cB$ for some small constant $c>0$. In addition $\langle a \rangle:=1+|a|$ and
$a\pm:=a\pm \epsilon$ with $0 < \epsilon <<1$.

\subsection{Norms and Strichartz estimates}
We use $L^r_x(\R^n)$ to denote the Lebesgue space of functions $f :
\R^n \rightarrow \C$ whose norm $$\|f\|_{L^r_x}:=\left(
\int_{\R^n}|f(x)|^r dx \right)^{\frac{1}{r}}$$ is finite, with the
usual modification in the case $r=\infty.$ We also define the
space-time spaces $L^q_tL^r_x$ by $$\|u\|_{L^q_tL^r_x} :=
\left(\int_J \|u\|_{L^r_x}^q dt\right)^{\frac{1}{q}}$$for any
space-time slab $J \times \R^n,$ with the usual modification when
either $q$ or $r$ are infinity. When $q=r$ we abbreviate
$L^q_tL^r_x$ by $L^q_{t,x}.$

\begin{defn} A pair of exponents $(q,r)$ is called admissible  in
    $\R^{n}$ if
$$ \frac{2}{q}+\frac{n}{r} = \frac{n}{2}, \ \ \ 2 \leq q,r \leq \infty, \ \ \ (q,r,n) \ne (2, \infty, 2).$$
\end{defn}
We recall the following Strichartz estimate \cite{gv}, \cite{kt}.
\begin{prop}Let $(q,r)$ and $(\tilde{q},\tilde{r})$ be
any two admissible pairs. Suppose that $u$ is a solution to
\begin{align}
&iu_{t}+ \Delta u -G(x,t)=0, \; \; x \in J \times \R^{n} ,\nonumber\\
&u(x,0)=u_{0}(x).\nonumber
\end{align} Then we have the estimate \begin{equation}\label{S}
\|u\|_{L^q_tL^r_x (J\times\R^n)} \lesssim
\|u_0\|_{L^2(\R^n)}+\|G\|_{L^{\tilde q^{\prime}}_tL^{\tilde r^{\prime}}_x
(J\times\R^n)}
\end{equation}
with the prime exponents denoting H\"older dual exponents.
\end{prop}

We now define the spatial Fourier transform on $\R^n$ by
$$\hat{f}(\xi) := \int{e^{-2\pi ix\cdot \xi} f(x)dx}.$$ We also
define the fractional differentiation operator $|\nabla|^\alpha$
for any real $\alpha$ by $$\widehat{|\nabla|^\alpha u}(\xi):=
|\xi|^\alpha \hat{u}(\xi)$$ and analogously
$$\widehat{\langle\nabla\rangle^\alpha u}(\xi):=
\langle\xi \rangle^\alpha \hat{u}(\xi).$$ We then define the
inhomogeneous Sobolev space $H^s$ and the homogeneous Sobolev
space $\dot{H}^s$ by $$\|u\|_{H^s} = \|\langle\nabla\rangle^s
u\|_{L^2_x}; \ \ \ \ \|u\|_{\dot{H}^s} =\|\nabla|^s u\|_{L^2_x}.$$

\subsection{Nonlinearity}
As in \cite{vz1,vz2} we use the notation $F(z)=|z|^p z$,
$p=\frac{4}{n},$ for  the function that defines the nonlinearity
in \eqref{ivp1}. We compute the derivatives
$$F_z(z)= \frac{p+2}{2}|z|^p, \ \ \ \textrm{and} \ \ \ F_{\overline{z}}(z)= \frac{p}{2}|z|^p \frac{z}{\overline{z}}.$$
We denote by $F'$ the vector $(F_z,F_{\overline{z}})$. Also we adopt
the notation $$w \cdot F'(z) = w F(z) +
\overline{w}F_{\overline{z}}(z).$$
In particular, the following chain rule is valid
$$\nabla F(u)= \nabla u \cdot F'(u).$$
Clearly
$F'(z)= O(|z|^p)$ and we can estimate the modulus of continuity of
$F'$ as follows
\begin{equation} \label{MC1}
|F'(z) - F'(w)| \lesssim
|z-w|^{\min\{1,p\}}(|z|+|w|)^{p-\min\{1,p\}},
\end{equation} for all $z,w \in \mathbb{C}.$ By the fundamental theorem of calculus we have that
$$F(z+w)-F(z)=\int_{0}^{1}w \cdot F^{\prime}(z+\theta w)d \theta$$
and thus the following estimate holds true
$$F(z+w) = F(z) + O(|w||z|^p)+ O(|w|^{p+1})$$for all complex
values $z,w$.

We notice that in the case $n=3,4$ the nonlinearity $F$ is in
$C^{1,1}(\mathbb{C})$, while in the case $n \geq 5$, $F \in
C^{1,4/n}(\mathbb{C})$, that is $F'$ is only H\"older
continuous. Hence, to estimate our nonlinearity, we will need the
following fractional chain rules\footnote{The reader should keep in
mind that these rules will be used  with  $G=F'$. }.

\begin{lem}\label{Lip}Suppose that $G \in C^{0,1}(\mathbb{C}), $ and $\alpha \in
(0,1)$. Then for
$\frac{1}{r_{1}}+\frac{1}{r_{2}}=\frac{1}{r}$, with $1< r\leq r_{2} < \infty $
and $1<r_{1}\leq \infty$ we have
$$\||\nabla|^\alpha G(u)\|_{L^r_x} \lesssim \|G'(u)\|_{L^{r_{1}}_x} \||\nabla|^\alpha u\|_{L^{r_{2}}_x}.$$
\end{lem}
The proof of this lemma when $1<r_{1}<\infty$  can be found in \cite{chw}
and when $r_{1} = \infty$  in \cite{kpvkdv}.

\begin{lem} \label{Holder}Suppose that $G \in C^{0,\alpha}(\mathbb{C}), \alpha \in
(0,1).$  Then, for every $0 < \sigma < \alpha, 1<r<\infty,$ and
$\sigma/\alpha < \rho < 1$ we have
$$\||\nabla|^\sigma G(u)\|_{L^r_x} \lesssim
\||u|^{\alpha - \frac{\sigma}{\rho}}\|_{L^{r_1}_x}
\||\nabla|^\rho u\|^{\frac{\sigma}{\rho}}_{L^{\frac{\sigma}{\rho}r_2}_x},$$
provided
$\frac{1}{r}=\frac{1}{r_1}+\frac{1}{r_2},$ and
$(1-\frac{\sigma}{\alpha \rho})r_1 > 1.$
\end{lem}
The proof of this lemma can be found in \cite{v:thesis}.

Also the following estimates can be found in \cite{vz2}. We notice that
for these estimates to hold, it suffices to require that $F \in
C^1(\mathbb{C}).$
\begin{lem} Let $1<r,r_1,r_2 < \infty$ be such that
$\frac{1}{r}=\frac{1}{r_1}+\frac{1}{r_2}.$ Then, for any $0 < \nu
< s $ we have
\begin{align}\label{diff}
  &\|\nabla I F(u) - (I \nabla u)F'(u)\|_{L^r_x} \lesssim N^{-1+s-\nu}\|\nabla I u\|_{L^{r_1}_x}
  \|\langle \nabla \rangle^{1-s+\nu}F'(u)\|_{L^{r_2}_x}, \\
  \label{oneterm}
 &\|\nabla I F(u)\|_{L^r_x} \lesssim  \|\nabla I u\|_{L^{r_1}_x} \|F'(u)\|_{L^{r_2}_x}+N^{-1+s-\nu}\|\nabla I u\|_{L^{r_1}_x}
  \|\langle \nabla \rangle^{1-s+\nu}F'(u)\|_{L^{r_2}_x}.
\end{align}
\end{lem}

\subsection{Littlewood-Paley Theory and the $I$-operator}

We shall also need some Littlewood-Paley theory. In particular, let
$\eta(\xi)$ be a smooth bump function supported in the ball $|\xi|
\leq 2,$ which is equal to one on the unit ball. Then, for each
dyadic number $M$ we define the Littlewood-Paley operators

\begin{align*}\widehat{P_{\leq M} f}(\xi) &= \eta(\xi/M) \hat{f}(\xi),\\
\widehat{P_{> M} f}(\xi) &= (1-\eta(\xi/M)) \hat{f}(\xi),\\
\widehat{P_{M} f}(\xi) &= (\eta(\xi/M)-\eta(2\xi/M)) \hat{f}(\xi).\\
\end{align*}
Similarly, we can define $P_{<M}, P_{\geq M}.$

Finally, we introduce the $I$-operator. For $s<1$ and a parameter $N
>>1$ let $m(\xi)$ be the following smooth monotone multiplier:
\[m(\xi):= \left\{\begin{array}{ll}
1 & \mbox{if $|\xi|<N$,}\\
(\frac{|\xi|}{N})^{s-1} & \mbox{if $|\xi|>2N$.}
\end{array}
\right.\] We define the multiplier operator $I:H^{s} \rightarrow
H^{1}$ by
$$\widehat{Iu}(\xi)=m(\xi)\hat{u}(\xi).$$ Some basic properties of
this operator are collected in the following Lemma.

\begin{lem}\label{multiplier}Let $1<r<\infty$ and $0<s<1$. Then,
\begin{align}
\label{mult1}& \|\langle \nabla\rangle^{\rho}P_{>N}f\|_{L^r_x} \lesssim
N^{\rho -1}\|\nabla I f\|_{L^r_x},\\
\label{mult2}& \|\langle \nabla\rangle^{\rho}f\|_{L^r_x} \lesssim
\|\langle \nabla \rangle I f\|_{L^r_x},\\\label{Ioperator}
&\|f\|_{H^s_x} \lesssim \|If\|_{H^1_x} \lesssim
N^{1-s}\|f\|_{H^s_x}.
\end{align}for all $0 \leq \rho \leq s.$
\end{lem}\begin{proof}We write,
$$ \|\langle \nabla \rangle^{\rho}P_{>N}f\|_{L^r_x} = \|P_{>N}\langle \nabla \rangle^{\rho}(\nabla I)^{-1}\nabla I
f\|_{L^r_x},$$and the claim \eqref{mult1} follows from H\"ormander's
multiplier theorem.

In order to get \eqref{mult2} we write
\begin{align}\|\langle
\nabla\rangle^{\rho}f\|_{L^r_x} \leq \|P_{\leq N}\langle
\nabla\rangle^{\rho}f\|_{L^r_x}+\|P_{>N}\langle
\nabla\rangle^{\rho}f\|_{L^r_x} \leq \| \langle \nabla\rangle
If\|_{L^r_x}+\|P_{>N}\langle \nabla\rangle^{\rho}f\|_{L^r_x}
\end{align}
and again the claim follows from H\"ormander's multiplier theorem
and \eqref{mult1}.
\\
\\
Finally, to show \eqref{Ioperator} we observe that by the
definition of the $I$-operator and \eqref{mult1} we get
\begin{align}
\|f\|_{H^s_x} &\lesssim \|P_{\leq N}f\|_{H^s_x} +
 \|\langle \nabla\rangle^s P_{>N}f\|_{L^2_x}\\\nonumber
& \lesssim \|I P_{\leq N}f\|_{H^s_x} + N^{s-1}\|\nabla I
f\|_{L^2_x}  \lesssim \|If\|_{H^1_x}.
\end{align}Furthermore,
\begin{equation}\|If\|_{H^1_x} = \|\langle \nabla \rangle^{1-s}I \langle \nabla \rangle^s
f\|_{L^2_x}\lesssim N^{1-s}\| \langle \nabla \rangle^s
f\|_{L^2_x}\lesssim N^{1-s}\|f\|_{H^s_x},
\end{equation}which concludes the proof of \eqref{Ioperator}.
\end{proof}

\subsection{Interaction Morawetz estimates}

We conclude this section with some interaction  Morawetz estimates. In \cite{ckstt6}
it was proved that when $n=3$, on any space-time slab $J\times \R^n$ on which the
solution $u$ to \eqref{ivp1}-\eqref{bc1} exists, the following
\emph{a priori} bound is satisfied:
\begin{equation} \label{Morawetz3}
\int_{J} \int_{\R^3} |u(x,t)|^4 dx dt \lesssim
\|u\|_{L^{\infty}_tL^2_x}^{2}\|u\|_{ L^{\infty}_t
\dot{H}^{\frac{1}{2}}_x (J \times \R^3) }^2 \lesssim
\|u_{0}\|_{L_{x}^{2}}^{2}\|u\|_{ L^{\infty}_t
\dot{H}^{\frac{1}{2}}_x (J \times \R^3) }^2.
\end{equation}
By H\"older's inequality in time we get
\begin{equation}\label{Mor0}
\|u\|_{L_{t}^{\frac{8}{3}}L_{x}^{4}} \lesssim
T^{\frac{1}{8}}\|u_{0}\|_{L_{x}^{2}}^{2}\|u\|_{ L^{\infty}_t
\dot{H}^{\frac{1}{2}}_x (J \times \R^3) }^2.
\end{equation}
A generalization of \eqref{Morawetz3}  in any dimension $n \geq 4$ was
proved in  \cite{tvz}:
\begin{align}\label{Morawetz}
\int_{J}\int_{\R^n}\int_{\R^n}\frac{|u(x,t)|^2|u(y,t)|^2}{|x-y|^3}dx
dy dt +
\int_{J}\int_{\R^n}\int_{\R^n}\frac{|u(x,t)|^2|u(y,t)|^{\frac{4}{n}+2}}{|x-y|}dx
dy dt & \\ \lesssim
\|u_0\|_{L^2_x}\|u\|^2_{L^\infty_t\dot{H}^{\frac{1}{2}}_x(J\times
\R^n ).}& \nonumber
\end{align}
As a consequence of \eqref{Morawetz} and some harmonic analysis
\cite{tvz},  gives the following \emph{a priori}  estimate for the solution to
\eqref{ivp1}-\eqref{bc1} for $n\geq 4$,
\begin{equation}\label{Mor2}
\||\nabla|^{-\frac{n-3}{4}}u\|_{L^4_{t,x}(J\times \R^n)} \lesssim
\|u_0\|_{L^2_x}^{\frac{1}{2}}\|u\|_{L^\infty_t\dot{H}^{\frac{1}{2}}(J\times
\R^n )}^{\frac{1}{2}}.
\end{equation}
Interpolating between \eqref{Mor2} and the trivial estimate
\begin{equation}
\||\nabla|^{\frac{1}{2}}u\|_{L_{t}^{\infty}L_{x}^{2}} \leq
\|u\|_{L^\infty_t \dot{H}^{\frac{1}{2}}_x,}
\end{equation}
we have that
\begin{equation}\label{Mor3}
\|u\|_{L^{2(n-1)}_tL^{\frac{2(n-1)}{n-2}}_x(J\times \R^n )}
\lesssim
\|u_0\|_{L^2_x}^{\frac{1}{2}}\|u\|_{L^\infty_t\dot{H}^{\frac{1}{2}}_x(J\times
\R^n )}^{\frac{n-2}{n-1}}.
\end{equation}
Finally, applying H\"older inequality in time, for $J=[0,T
]$, we obtain
\begin{equation}\label{Morfinal}
\|u\|_{L^{\frac{4(n-1)}{n}}_tL^{\frac{2(n-1)}{n-2}}_x(J\times \R^n
)} \lesssim T^{\frac{n-2}{4(n-1)}}
\|u_0\|_{L^2_x}^{\frac{1}{2}}\|u\|_{L^\infty_t\dot{H}^{\frac{1}{2}}_x(J\times
\R^n )}^{\frac{n-2}{n-1}}.
\end{equation}
Notice that \eqref{Mor0} coincides with the inequality obtained by
formally substituting $n=3$ into \eqref{Morfinal}. Thus from now on we may use
\eqref{Morfinal}, for every $n \geq 3$.

Also we remark that the pair $(\frac{4(n-1)}{n},\frac{2(n-1)}{n-2})$ is
admissible.

\section{The main estimates and the local well-posedness for the $I$-system.}

Assume $J=[a,b]$. We denote by \begin{equation}\label{ZI}
Z_I(J) :=
\sup_{(q,r) \ \ admissible} \|\langle \nabla \rangle I
u\|_{L^q_tL^r_x(J\times \R^n)}.
\end{equation}
Often we will drop the dependence on the time interval $J$, and we
will write $Z_I.$

Our main estimates for $F(u)=|u|^{\frac{4}{n}}u$ read as
follow\footnote{ Similar types of estimates are also proved in
\cite{vz1, vz2} and in a weaker version  in \cite{ckstt7}.}

\begin{prop}\label{3.1} Let $1>s>\frac{1}{1+\min\{1,\frac{4}{n}\}}$, and let $0 < \nu \leq \min\{1,\frac{4}{n}\}s+s -1.$ Then,
\begin{align}\label{diff1decay}
 &\|\nabla I F(u) - (I \nabla u)F'(u)\|_{L^{\frac{4(n-1)}{3n-4}}_tL^{\frac{2(n-1)}{n}}_x}
\lesssim N^{-1+s-\nu}Z_I^{1+\frac{4}{n}},\\\label{one1decay}
&\|\nabla I
F(u)\|_{L^{\frac{4(n-1)}{3n-4}}_tL^{\frac{2(n-1)}{n}}_x} \lesssim
\|u\|^{\frac{4}{n}}_{L^{\frac{4(n-1)}{n}}_tL^{\frac{2(n-1)}{n-2}}_x}Z_I+
N^{-1+s-\nu}Z_I^{1+\frac{4}{n}},
\\\label{diff2decay}
 &\|\nabla I F(u) - (I \nabla u)F'(u)\|_{L^{2}_tL^{\frac{2n}{n+2}}_x}
\lesssim N^{-1+s-\nu}Z_I^{1+\frac{4}{n}},\\
\label{one2bound}
   &\|\nabla I F(u)\|_{L^{2}_tL^{\frac{2n}{n+2}}_x} \lesssim
   Z_I^{1+\frac{4}{n}}.\\\nonumber
\end{align}
\end{prop}

\begin{proof}
First we prove the estimate \eqref{diff1decay}. In order to do that we apply
\eqref{diff} to the left hand side of \eqref{diff1decay} with
$r=2(n-1)/n$, $r_1=2(n-1)/(n-2)$ and $r_2=n-1$.
Performing H\"older's inequality in time, we then obtain
\begin{align}\label{A1}
&\|\nabla I F(u) - (I \nabla
u)F'(u)\|_{L^{\frac{4(n-1)}{3n-4}}_tL^{\frac{2(n-1)}{n}}_x}
\\\nonumber &\lesssim  N^{-1+s-\nu} \|\nabla I
u\|_{L^{\frac{4(n-1)}{n}}_tL^{\frac{2(n-1)}{n-2}}_x} \|\langle
\nabla \rangle^{1-s+\nu}
F'(u)\|_{L^{\frac{2(n-1)}{n-2}}_tL^{n-1}_x}.
\end{align}
Notice that the pair $(\frac{4(n-1)}{n},\frac{2(n-1)}{n-2})$ is
admissible, hence
\begin{equation}\label{A2}
\|\nabla I F(u) - (I \nabla
u)F'(u)\|_{L^{\frac{4(n-1)}{3n-4}}_tL^{\frac{2(n-1)}{n}}_x}
\lesssim  N^{-1+s-\nu}Z_I \|\langle \nabla \rangle^{1-s+\nu}
F'(u)\|_{L^{\frac{2(n-1)}{n-2}}_tL^{n-1}_x}.
\end{equation}Thus, we need to control $\|\langle \nabla \rangle^{1-s+\nu}
F'(u)\|_{L^{\frac{2(n-1)}{n-2}}_tL^{n-1}_x}.$ Towards this aim, we
need to distinguish the cases $n=3, 4$ and the case $n\geq 5.$
\\
\\
\noindent\textbf{Case $\mathbf{n=3}$.} We  first write
\begin{align}\label{n3}
&\|\langle \nabla \rangle^{1-s+\nu}F'(u)\|_{L^{4}_tL^{2}_x} \lesssim
\||\nabla|^{1-s+\nu}F'(u)\|_{L^{4}_tL^{2}_x}+
\|F'(u)\|_{L^{4}_tL^{2}_x}.
\end{align}
Since $F'(u)=O(|u|^{\frac{4}{3}})$ we have
\begin{equation}\label{n3a}
\|F'(u)\|_{L^{4}_tL^{2}_x}
\lesssim \|u\|_{L^{\frac{16}{3}}_tL^{\frac{8}{3}}_x}^{\frac{4}{3}}
\lesssim Z_{I}^{\frac{4}{3}},
\end{equation}
where the last inequality follows from the fact that the pair
$(\frac{16}{3},\frac{8}{3})$ is admissible together with Lemma
2.6. To bound the first term on the right-hand side of \eqref{n3}
we use Lemma \ref{Lip} and obtain for $1/q_{1}+1/q_{2}=1/4$ and
$1/r_{1}+1/r_{2}=1/2$,
\begin{equation} \label{n3der}
\||\nabla|^{1-s+\nu}F'(u)\|_{L^{4}_tL^{2}_x}\lesssim
\|u\|_{L^{\frac{q_{1}}{3}}_tL^{\frac{r_{1}}{3}}_x}^{1/3}
\||\nabla|^{1-s+\nu}u\|_{L^{q_{2}}_tL^{r_{2}}_x}.
\end{equation}
Now if the pairs $(\frac{q_{1}}{3},\frac{r_{1}}{3})$ and
$(q_{2},r_{2})$ are admissible then by Lemma \ref{multiplier}
the expression \eqref{n3der} implies
\begin{equation} \label{n3b}
\||\nabla|^{1-s+\nu}F'(u)\|_{L^{4}_tL^{2}_x}\lesssim Z_{I}^{\frac{4}{3}},
\end{equation}
as long as $0 \leq \nu \leq 2s -1$. This is guaranteed by the
assumption $0 < \nu \leq \min\{1,\frac{4}{n}\}s+s -1.$ The above
pairs are admissible if
$$\frac{2}{q_{2}}+\frac{3}{r_{2}}=\frac{3}{2}\ \ \ ,\
\frac{2}{q_{1}}+\frac{3}{r_{1}}=\frac{1}{2}.$$
If we add these two equalities we get
$$\frac{2}{q_{1}}+\frac{2}{q_{2}}+\frac{3}{r_{1}}+\frac{3}{r_{2}}=2,$$
which is exactly the condition on the H\"older's exponents. Thus
by combining \eqref{A2} with \eqref{n3}, \eqref{n3a} and
\eqref{n3b}, we obtain the desired estimate \eqref{one1decay}.
\\
\\
\noindent\textbf{Case $\mathbf{n=4}$.} Now $F'(u)=O(|u|)$, so
\begin{equation}\label{n4}
\|\langle \nabla \rangle^{1-s+\nu}F'(u)\|_{L^{3}_tL^{3}_x}
\lesssim \||\nabla|^{1-s+\nu}F'(u)\|_{L^{3}_tL^{3}_x}+
\|F'(u)\|_{L^{3}_tL^{3}_x}.
\end{equation}
Since the pair $(3,3)$ is admissible, again using Lemma 2.6 we
have that
\begin{equation} \label{n4u}
\|F'(u)\|_{L^{3}_tL^{3}_x}\lesssim \|u\|_{L^{3}_tL^{3}_x} \lesssim Z_{I}.
\end{equation}
To bound the first term on the right-hand side of \eqref{n4} we
use Lemma \ref{Lip} with $r,r_2=3$, and $r_1=\infty$. Hence,
\begin{equation} \label{n4nab}
\||\nabla|^{1-s+\nu}F'(u)\|_{L^{3}_tL^{3}_x}
\lesssim \||\nabla|^{1-s+\nu}u\|_{L^{3}_tL^{3}_x}
\leq Z_I,
\end{equation}
where the last inequality follows from Lemma \ref{multiplier}, as
long as $0 \leq \nu \leq 2s -1$. This is guaranteed by the
assumption $0 < \nu \leq \min\{1,\frac{4}{n}\}s+s -1.$ Thus by
combining \eqref{A2} with \eqref{n4}, \eqref{n4u} and
\eqref{n4nab} we obtain the desired estimate \eqref{one1decay}.
\\
\\
\noindent\textbf{Case $\mathbf{n\geq 5}$.} First we bound
$\|F'(u)\|_{L^{\frac{2(n-1)}{n-2}}_tL^{n-1}_x}$, where
$F'(u)=O(|u|^{\frac{4}{n}}$). Hence
$$\|F'(u)\|_{L^{\frac{2(n-1)}{n-2}}_tL^{n-1}_x}\lesssim
\|u\|^{4/n}_{L^{\frac{8(n-1)}{n(n-2)}}_tL^{\frac{4(n-1)}{n}}_x},$$
where the pair $(\frac{8(n-1)}{n(n-2)},\frac{4(n-1)}{n})$ is
admissible. Thus  we immediately have that
\begin{equation} \label{n5inh}
\|F'(u)\|_{L^{\frac{2(n-1)}{n-2}}_tL^{n-1}_x}\lesssim Z_{I}^{\frac{4}{n}}.
\end{equation}
To bound the homogeneous derivative, since $4/n<1$,  we apply Lemma \ref{Holder},
with $\alpha=4/n, \sigma=1-s+\nu, r=n-1$ and $r_1,r_2$ satisfying
$$(\frac{4}{n} - \frac{\sigma}{\rho})r_1 = \frac{4(n-1)}{n}, \, \, \,
\frac{\sigma r_2}{\rho} = \frac{4(n-1)}{n},$$ with $\sigma
\frac{4}{n} < \rho <1$ to be chosen later. Notice that in order to
apply such lemma, we need $\left(1-\frac{\sigma}{\alpha
\rho}\right)r_{1}>1$. For our choices of values, this quantity
equals $n-1$, hence the required assumption is satisfied. Then,
$$\||\nabla|^{1-s+\nu}
F'(u)\|_{L^{n-1}_x} \lesssim
\|u\|^{\epsilon_1}_{L^{\frac{4(n-1)}{n}}_x} \||\nabla|^{\rho}u
\|^{\epsilon_2}_{L^{\frac{4(n-1)}{n}}_x},$$
where
$$ \epsilon_1 = \frac{4}{n} - \frac{\sigma}{\rho}, \, \, \, \epsilon_2=\frac{\sigma}{\rho},  \ \
\text{hence} \ \  \epsilon_1+\epsilon_2 =\frac{4}{n}.$$ Thus,
applying H\"older's inequality in time, we obtain
\begin{equation*}
\||\nabla|^{1-s+\nu} F'(u)\|_{L^{\frac{2(n-1)}{n-2}}_tL^{n-1}_x}
\lesssim
\|u\|^{\epsilon_1}_{L^{q_1\epsilon_1}_tL^{\frac{4(n-1)}{n}}_x}
\||\nabla|^{\rho}u
\|^{\epsilon_2}_{L^{q_2\epsilon_2}_tL^{\frac{4(n-1)}{n}}_x},
\end{equation*}
where $$q_i = \frac{8(n-1)}{n(n-2)} \frac{1}{\epsilon_i}, \, \,
i=1,2.$$
Finally, for any $\rho$ such that $(1-s+\nu)\frac{4}{n}
<\rho < 1$, which exists by our assumptions on $\nu$ since $\nu \leq
\frac{4}{n}s+s-1$, we have
\begin{align}
\||\nabla|^{1-s+\nu} F'(u)\|_{L^{\frac{2(n-1)}{n-2}}_tL^{n-1}_x}
& \lesssim \|u\|^{\frac{4}{n} - \frac{\sigma}{\rho}}_
{L^{\frac{8(n-1)}{n(n-2)}}_tL^{\frac{4(n-1)}{n}}_x}
\||\nabla|^{\rho}u\|^{\frac{\sigma}{\rho}}_{L^{\frac{8(n-1)}{n(n-2)}}_tL^{\frac{4(n-1)}{n}}_x}
\nonumber \\
& \lesssim \|u\|^{\frac{4}{n} - \frac{\sigma}{\rho}}_
{L^{\frac{8(n-1)}{n(n-2)}}_tL^{\frac{4(n-1)}{n}}_x}
\|\langle \nabla \rangle ^{\rho}u\|^{\frac{\sigma}{\rho}}_
{L^{\frac{8(n-1)}{n(n-2)}}_tL^{\frac{4(n-1)}{n}}_x}. \label{B1}
\end{align}
The pair $(\frac{8(n-1)}{n(n-2)},\frac{4(n-1)}{n})$ is admissible.
Hence,
the expression \eqref{B1} implies by Lemma \ref{multiplier} that
\begin{equation}\label{B2}
\|\langle \nabla \rangle^{1-s+\nu}
F'(u)\|_{L^{\frac{2(n-1)}{n-2}}_tL^{n-1}_x} \lesssim
Z_I^{\frac{4}{n}},
\end{equation} as long as $0 \leq \rho \leq s.$ By our assumption on $\nu$, there exists $\rho$
such that $(1-s+\nu)\frac{4}{n} <\rho \leq s$. Then, combining
\eqref{A2}, \eqref{n5inh} and \eqref{B2}, we obtain
$$
\|\nabla I F(u) - (I \nabla
u)F'(u)\|_{L^{\frac{4(n-1)}{3n-4}}_tL^{\frac{2(n-1)}{n}}_x}
\lesssim  N^{-1+s-\nu}Z_I ^{1+\frac{4}{n}},
$$
which is the desired estimate.

We now proceed with the proof of \eqref{one1decay}. By the triangle
inequality, and the estimate \eqref{diff1decay}, we obtain
\begin{align*}
& \|\nabla IF(u)\|_{L^{\frac{4(n-1)}{3n-4}}_tL^{\frac{2(n-1)}{n}}_x} \\
& \lesssim \|(I\nabla u)F'(u)\|_{L^{\frac{4(n-1)}{3n-4}}_tL^{\frac{2(n-1)}{n}}_x}
+ \|\nabla I F(u) - (I \nabla u)F'(u)\|_{L^{\frac{4(n-1)}{3n-4}}_tL^{\frac{2(n-1)}{n}}_x} \\
& \lesssim \|(I\nabla u) F'(u)\|_{L^{\frac{4(n-1)}{3n-4}}_tL^{\frac{2(n-1)}{n}}_x}
+ N^{-1+s-\nu}Z_I ^{1+\frac{4}{n}}.
\end{align*}
In order to conclude the proof of \eqref{one1decay}, we need to estimate $\|(I\nabla u)
F'(u)\|_{L^{\frac{4(n-1)}{3n-4}}_tL^{\frac{2(n-1)}{n}}_x}.$
H\"older's inequality gives,
\begin{align}\nonumber \|(I\nabla u)
F'(u)\|_{L^{\frac{4(n-1)}{3n-4}}_tL^{\frac{2(n-1)}{n}}_x}
&\lesssim
\||u|^{\frac{4}{n}}\|_{L^{n-1}_tL^{\frac{(n-1)n}{2(n-2)}}_x}\|I\nabla
u\|_{L^{\frac{4(n-1)}{3n-8}}_tL^{\frac{2n(n-1)}{n^2-4n+8}}_x}\\&\lesssim
\|u\|^{\frac{4}{n}}_{L^{\frac{4(n-1)}{n}}_tL^{\frac{2(n-1)}{n-2}}_x}Z_I,
\end{align}where the last inequality follows from the fact that
the pair $(\frac{4(n-1)}{3n-8},\frac{2n(n-1)}{n^2-4n+8})$ is
admissible. This concludes the proof of \eqref{one1decay}.

The proof of \eqref{diff2decay} is along the lines of
\eqref{diff1decay}. Indeed, by \eqref{diff} we
get,
\begin{align}\label{A4}
& \|\nabla I F(u) - (I \nabla
u)F'(u)\|_{L^{2}_tL^{\frac{2n}{n+2}}_x}
\\ \nonumber
& \lesssim  N^{-1+s-\nu} \|\nabla I u\|_{L^{q_1}_tL^{r_1}_x}
\|\langle \nabla \rangle^{1-s+\nu} F'(u)\|_{L^{q_2}_tL^{r_2}_x},
\end{align}
where
\begin{equation}\label{q2r2}\frac{1}{2}=\frac{1}{q_1}+\frac{1}{q_2}, \frac{n+2}{2n}=
\frac{1}{r_1}+\frac{1}{r_2},  \ \ \text{and} \ \
\frac{2}{q_2}+\frac{n}{r_2}=2.
\end{equation} The last condition follows from choosing the pair
$(q_1,r_1)$ to be admissible. Hence, setting $$r_2=n-1, \ \
\text{and} \ \ q_2=\frac{2(n-1)}{n-2},$$ we obtain
\begin{equation}\label{A5}
 \|\nabla I F(u) - (I \nabla
u)F'(u)\|_{L^{2}_tL^{\frac{2n}{n+2}}_x}
 \lesssim  N^{-1+s-\nu} Z_I \|\langle \nabla \rangle^{1-s+\nu}
F'(u)\|_{L^{\frac{2(n-1)}{n-2}}_tL^{n-1}_x}.
\end{equation}Notice that the right-hand side of \eqref{A5} was
estimated in the proof of \eqref{diff1decay}. This concludes the
proof of \eqref{diff2decay}.

We now turn to \eqref{one2bound}. The triangle inequality and
\eqref{diff2decay} imply,

\begin{equation}
\|\nabla I F(u)\|_{L^{2}_tL^{\frac{2n}{n+2}}_x}
 \lesssim  \|(\nabla I u)F'(u)\|_{L^{2}_tL^{\frac{2n}{n+2}}_x} + N^{-1+s-\nu} Z_I^{1+\frac{4}{n}}.
\end{equation} Let $(q_1,r_1)$ be an admissible pair, and let
$(q_2,r_2)$ be as in \eqref{q2r2}. Then H\"older's inequality
yields,
\begin{equation}
\|(\nabla I u)F'(u)\|_{L^{2}_tL^{\frac{2n}{n+2}}_x} \lesssim
\|\nabla I u\|_{L^{q_1}_tL^{r_1}_x}\|F'(u)\|_{L^{q_2}_tL^{r_2}_x}
\lesssim Z_I
\|u\|^{\frac{4}{n}}_{L^{\frac{4q_2}{n}}_tL^{\frac{4r_2}{n}}_x}\lesssim
Z_I^{1+\frac{4}{n}},
\end{equation} where the last inequality follows from the fact
that $(\frac{4q_2}{n},\frac{4r_2}{n})$ is admissible.
\end{proof}

\begin{rem}We notice that the estimates in Proposition \ref{3.1}
hold true for any choice of dual admissible pair $(q',r')$ on the
right-hand side. This follows immediately by a simple modification
of the proof above.
\end{rem}

We now turn to the proof of a ``modified" local existence theory,
that is a local existence involving norms of $Iu$ instead of $u$.
Following for example the argument in \cite{FG}, the proof can be
reduced to showing the following:
\begin{prop}\label{lwp}
Let $1>s>\frac{1}{1+\min\{1,\frac{4}{n}\}}$ and assume that if $u$ is a
solution to \eqref{ivp1}-\eqref{bc1} on $J=[a,b]$, the a priori
bound
$$\|u\|^{\frac{4(n-1)}{n}}_{L^{\frac{4(n-1)}{n}}_tL^{\frac{2(n-1)}{n-2}}_x(J \times \R^n)}
< \mu,$$
\\
with $\mu$ a small universal constant, is satisfied. Then $Z_I(J) \lesssim
\|Iu(a)\|_{H^1}$.
\end{prop}
\begin{proof}
Applying $I\langle \nabla \rangle$ to \eqref{ivp1}, and using the
Strichartz estimate in \eqref{S}, for any pair of admissible
exponents $(q,r)$ we have
\begin{align}
\|\langle \nabla \rangle I u\|_{L^q_tL^r_x}
& \lesssim \|Iu(a)\|_{H^1} + \|\langle \nabla  \rangle
I F(u)\|_{L^{\frac{4(n-1)}{3n-4}}_tL^{\frac{2(n-1)}{n}}_x} \nonumber \\
& \lesssim \|Iu(a)\|_{H^1} +
\|IF(u)\|_{L^{\frac{4(n-1)}{3n-4}}_tL^{\frac{2(n-1)}{n}}_x}+N^{-1+s+\nu}Z_I^{1+\frac{4}{n}}
+ \mu^{\frac{1}{n-1}}Z_I, \label{lwp2}
\end{align}
where the last inequality follows from \eqref{one1decay}, and $0 <
\nu \leq \min \{1,\frac{4}{n}\}s+s-1.$

We now need to control
$\|IF(u)\|_{L^{\frac{4(n-1)}{3n-4}}_tL^{\frac{2(n-1)}{n}}_x}.$ We
compute,
\begin{equation}
\|IF(u)\|_{L^{\frac{4(n-1)}{3n-4}}_tL^{\frac{2(n-1)}{n}}_x} \leq
\|IP_{<N}F(u)\|_{L^{\frac{4(n-1)}{3n-4}}_tL^{\frac{2(n-1)}{n}}_x}
+\|IP_{\geq
N}F(u)\|_{L^{\frac{4(n-1)}{3n-4}}_tL^{\frac{2(n-1)}{n}}_x}.
\end{equation}
Using standard harmonic analysis we can bound

\begin{equation*}\|IP_{\geq
N}F(u)\|_{L^{\frac{4(n-1)}{3n-4}}_tL^{\frac{2(n-1)}{n}}_x}\leq
\frac{1}{N}\|\nabla
IF(u)\|_{L^{\frac{4(n-1)}{3n-4}}_tL^{\frac{2(n-1)}{n}}_x},\end{equation*}
which in turn is bounded in \eqref{one1decay}. Moreover, by
H\"older's inequality we have,
\begin{align*}
\|IP_{<N}F(u)\|_{L^{\frac{4(n-1)}{3n-4}}_tL^{\frac{2(n-1)}{n}}_x}
& =\|F(u)\|_{L^{\frac{4(n-1)}{3n-4}}_tL^{\frac{2(n-1)}{n}}_x} \\
& \lesssim
\|u\|^{\frac{4}{n}}_{L^{\frac{4(n-1)}{n}}_tL^{\frac{2(n-1)}{n-2}}_x}
\|u\|_{L^{\frac{4(n-1)}{3n-8}}_tL^{\frac{2n(n-1)}{n^2-4n+8}}_x}
\lesssim \mu^{\frac{1}{n-1}}Z_I,
\end{align*}
where the last inequality follows from Lemma \ref{multiplier} together with the
fact that the pair
$(\frac{4(n-1)}{3n-8},\frac{2n(n-1)}{n^2-4n+8})$ is admissible.
Thus,
\begin{equation}
\|IF(u)\|_{L^{\frac{4(n-1)}{3n-4}}_tL^{\frac{2(n-1)}{n}}_x}
\lesssim N^{-1+s+\nu}Z_I^{1+\frac{4}{n}} + \mu^{\frac{1}{n-1}}Z_I,
\end{equation}
which combined with \eqref{lwp2} gives
\begin{equation}\label{lwp3}
Z_I \lesssim \|Iu(a)\|_{H^1} + N^{-1+s+\nu}Z_I^{1+\frac{4}{n}} +
\mu^{\frac{1}{n-1}}Z_I.
\end{equation}
A standard continuity argument finishes the proof if we pick $\mu$ sufficiently small and $N$ sufficiently
 large.
\end{proof}

\section{Almost conservation of the modified energy.}

In this section we will prove that the modified energy functional
$E(Iu)$ is almost conserved. We recall that,
$$E(u)(t)=\frac{1}{2}\int |\nabla u(t)|^{2}dx+\frac{n}{2n+4}\int |u(t)|^{\frac{2n+4}{n}}dx=E(u)(0)$$ for
any smooth solution $u$ to \eqref{ivp1}.We wish to prove the
following decay for the increment of the modified energy\footnote{
A similar estimate was also proved in \cite{vz1,vz2}.}.

\begin{prop} \label{endecay}Let $1>s>\frac{1}{1+\min\{1,\frac{4}{n}\}}$
and let $0 < \nu \leq \min \{1,\frac{4}{n}\}s+s-1.$
Assume that $u\in C_0^{\infty}(\R^n)$ is a solution to
\eqref{ivp1} on a time interval $[0,T_0]$. Then, for any
$J=[a,b]\subset [0,T_0]$
\begin{equation}\label{decay}
|E(Iu)(b) - E(Iu)(a)| \lesssim N^{-1+s-\nu} Z_I^{2+\frac{8}{n}}(J).
\end{equation}
\end{prop}

\begin{proof} Without loss of generality let us
assume that $J=[0,T]$. By the Fundamental Theorem of Calculus
\begin{align}
E(Iu)(T) - E(Iu)(0) & = \int_{0}^{T}\frac{\partial}{\partial t}E(Iu)(t)\, dt\\
&= Re\int_{0}^{T}\int_{\R^n}\overline{Iu_t}(-\Delta Iu +F(Iu))\, dx \, dt.
\end{align}
Since $Iu_t = i\Delta Iu - iIF(u)$, we get
$$Re\int_{0}^{T}\int_{\R^n}\overline{Iu_t}(-\Delta Iu +I F(u))\, dx \, dt=0.$$
Hence, after integration by parts we obtain
\begin{align}\nonumber
E(Iu)(T) - E(Iu)(0)
& = - \textrm{Im} \int_{0}^{T}\int_{\R^n}
\overline{\nabla Iu} \cdot \nabla \left[F(Iu)-IF(u)\right] \, dx \, dt\\ \nonumber
& - \textrm{Im} \int_{0}^{T}\int_{\R^n}
\overline{IF(u)} \cdot \left[F(Iu)-IF(u) \right] \, dx \, dt\\\label{energy1}
& = - \textrm{Im}\int_{0}^{T}\int_{\R^n}
\overline{\nabla Iu} \cdot (\nabla Iu) \left[F'(Iu)-F'(u) \right]\, dx \, dt\\\label{energy2}
& - \textrm{Im}\int_{0}^{T}\int_{\R^n}
\overline{\nabla Iu} \cdot \left[(\nabla Iu)F'(u) - \nabla I F(u) \right] \, dx \, dt\\\label{energy3}
& -\textrm{Im}
\int_{0}^{T}\int_{\R^n}
\overline{IF(u)} \cdot \left[ F(Iu)-IF(u) \right] \, dx \, dt.
\end{align}

We start by estimating \eqref{energy1}. By H\"older's inequality we get,
\begin{align}
|\eqref{energy1}|
& \lesssim \|\nabla I u\|_{L^2_t L^{\frac{2n}{n-2}}_x}
\| (\nabla Iu) \; \left[F'(Iu)-F'(u)\right]\|_{L^{2}_t L^{\frac{2n}{n+2}}_x} \nonumber \\
& \lesssim Z_I \;
\| (\nabla Iu) \; \left[F'(Iu)-F'(u)\right]\|_{L^{2}_t L^{\frac{2n}{n+2}}_x} \label{inctr1},
\end{align}
where to obtain \eqref{inctr1} we use the fact that the pair
$(2, \frac{2n}{n-2})$ is admissible. Now we proceed to bound
$\| (\nabla Iu) \; \left[F'(Iu)-F'(u)\right]\|_{L^{2}_t L^{\frac{2n}{n+2}}_x}$.
By H\"{o}lder's inequality we have
\begin{align}
\| (\nabla Iu) \; \left[F'(Iu)-F'(u)\right]\|_{L^{2}_t L^{\frac{2n}{n+2}}_x}
& \leq
\|\nabla I u\|_{L^2_t L^{\frac{2n}{n-2}}_x}
\|F'(Iu)-F'(u)\|_{L^{\infty}_t L^{\frac{n}{2}}_x} \nonumber \\
& \lesssim Z_I \;
\| \; |Iu - u|^{ \min\{1,\frac{4}{n}\} } \;
\left( |Iu| + |u| \right)^{\frac{4}{n} - \min\{1, \frac{4}{n}\} }
\|_{L^{\infty}_t L^{\frac{n}{2}}_x} \label{mvabs} \\
& \lesssim Z_I \; \|P_{>N}u\|^{\min\{1,\frac{4}{n}\}}_{L^{\infty}_t L^2_x}
\|u\|^{\frac{4}{n} - \min\{1, \frac{4}{n} \}}_{L^{\infty}_t L^2_x}, \label{mincases}
\end{align}
where to obtain \eqref{mvabs} we use the fact that the pair
$(2, \frac{2n}{n-2})$ is admissible and \eqref{MC1}, 
while to obtain \eqref{mincases} we distinguish the cases 
depending on $\min\{1, \frac{4}{n}\}$ and then use H\"{o}lder's inequality
in case $n \leq 4$. 
However \eqref{mult1} implies
$$ \|P_{>N} u \|_{L^{\infty}_t L^2_x} \lesssim N^{-1} \|\nabla I u\|_{L^{\infty}_t L^2_x}.$$
Thus
\begin{align}
\|P_{>N}u\|^{\min\{1,\frac{4}{n}\}}_{L^{\infty}_t L^2_x}
& \lesssim N^{-\min\{1, \frac{4}{n}\}}
\|\nabla I u\|_{L^{\infty}_t L^2_x}^{\min\{1, \frac{4}{n}\}} \nonumber \\
& \lesssim  N^{-\min\{1, \frac{4}{n}\}} Z_I^{\min\{1,
\frac{4}{n}\}} \label{en1hi} ,
\end{align}
where \eqref{en1hi} follows from the fact that the pair $(\infty, 2)$ is admissible.
On the other hand by splitting $u$ into high and low frequencies we obtain
\begin{align}
\|u\|_{L^{\infty}_t L^2_x}
& \leq \|P_{<N}u \|_{L^{\infty}_t L^2_x} +  \|P_{>N}u \|_{L^{\infty}_t L^2_x} \nonumber \\
& \lesssim \|I P_{<N}u \|_{L^{\infty}_t L^2_x} + N^{-1} \|\nabla Iu\|_{L^{\infty}_t L^2_x}
\label{en1uhl} \\
& \lesssim Z_I, \label{en1u}
\end{align}
where to obtain \eqref{en1uhl} we used the definition of the operator $I$ and \eqref{mult1},
while to obtain \eqref{en1u} we used the fact that the pair $(\infty, 2)$ is admissible.
However \eqref{mincases}, \eqref{en1hi} and \eqref{en1u} imply that
\begin{equation} \label{inctr1mvt}
\| (\nabla Iu) \; \left[F'(Iu)-F'(u)\right]\|_{L^{2}_t L^{\frac{2n}{n+2}}_x}
\lesssim N^{-\min\{1, \frac{4}{n}\}} Z_I^{1 + \frac{4}{n}}.
\end{equation}
Now we combine \eqref{inctr1} and \eqref{inctr1mvt} to conclude
$$ |\eqref{energy1}| \lesssim  N^{-\min\{1, \frac{4}{n}\}} Z_I^{2+ \frac{4}{n}},$$
which for  $1>s$ and $0 < \nu \leq \min \{1,\frac{4}{n}\}s+s-1$
implies
$$ |\eqref{energy1}| \lesssim N^{-1+s -\nu} Z_I^{2+ \frac{4}{n}}.$$

We now proceed to bound \eqref{energy2}. Applying H\"older's inequality we get
\begin{align}\nonumber
|\eqref{energy2}|
&\lesssim \|\nabla I u\|_{L^{\frac{4(n-1)}{n}}_tL^{\frac{2(n-1)}{n-2}}_x}
\|(\nabla Iu)F'(u)- \nabla I F(u) \|_{L^{\frac{4(n-1)}{3n-4}}_tL^{\frac{2(n-1)}{n}}_x}
\\\label{energy21}
& \lesssim Z_I
\|(I\nabla u)F'(u) - \nabla IF(u)\|_{L^{\frac{4(n-1)}{3n-4}}_tL^{\frac{2(n-1)}{n}}_x}\\\label{energy22}
&\lesssim N^{-1+s-\nu}Z_I^{2+\frac{4}{n}},
\end{align}
where to obtain \eqref{energy21} we used that the pair
$(\frac{4(n-1)}{n},\frac{2(n-1)}{n-2})$ is admissible, while to
obtain \eqref{energy22} we used \eqref{diff1decay}.

Finally, we bound \eqref{energy3}.
\begin{align}\nonumber
|\eqref{energy3}|
&\lesssim \||\nabla|^{-1}I(F(u))\|_{L^{2}_tL^{\frac{2n}{n-2}}_x}
\|\nabla \left[F(Iu)-IF(u)\right]\|_{L^{2}_tL^{\frac{2n}{n+2}}_x}\\\label{energy31}
& \lesssim \||\nabla I(F(u))\|_{L^{2}_tL^{\frac{2n}{n+2}}_x}
\|\nabla \left[F(Iu)-IF(u)\right]\|_{L^{2}_tL^{\frac{2n}{n+2}}_x}\\\label{energy32}
& \lesssim Z_I^{1+\frac{4}{n}}
\left(\|(\nabla Iu) \left[F'(Iu)-F'(u)\right]\|_{L^{2}_tL^{\frac{2n}{n+2}}_x}
+ \|(\nabla Iu)F'(u)- \nabla IF(u) \|_{L^{2}_tL^{\frac{2n}{n+2}}_x}\right)\\\label{energy33}
& = Z_I^{1+\frac{4}{n}}(
 N^{-\min\{1, \frac{4}{n}\}} Z_I^{1 + \frac{4}{n}}
+ N^{-1+s-\nu} {Z_I^{1+\frac{4}{n}}}),
\end{align}
where in order to obtain \eqref{energy31} we used Sobolev's
embedding Theorem, while to obtain \eqref{energy32} we used
\eqref{one2bound} and the triangle inequality, and to obtain
\eqref{energy33} we used \eqref{inctr1mvt} and \eqref{diff2decay}.
Finally, for  $1>s$ and $0 < \nu \leq \min \{1,\frac{4}{n}\}s+s-1$
\eqref{energy33} implies
$$ |\eqref{energy3}| \lesssim N^{-1+s -\nu} Z_I^{2+
\frac{8}{n}},$$which concludes our proof.
\end{proof}

\section{Global-well posedness.}

In this section we present the
proof of Theorem \ref{main}. We recall once again that this proof was
inspired by the arguments in \cite{FG}.
\begin{proof}We start by observing that the assumption on $s$
guarantees that $s > \frac{1}{1+\min\{1, \frac{4}{n}\}}$, thus we can apply the results
in the previous sections.

Now, suppose that $u(t,x)$ is a global in time solution to
\eqref{ivp1} with  initial data $u_0 \in C_0^\infty(\R^n)$. Set
$u^\lambda(x)=\frac{1}{\lambda^{\frac{n}{2}}}u(\frac{x}{\lambda},\frac{t}{\lambda^{2}})$.
We choose the parameter $\lambda$ so that $\|I u_0^\lambda\|_{H^1}
= O(1)$, that is
\begin{equation}\label{L}
\lambda \sim N^{\frac{1-s}{s}}.
\end{equation}
Next, let us pick a time
$T_0$ arbitrarily large, and  inspired by  \eqref{Morfinal} let us define
\begin{equation}
S : = \{0 < t < \lambda^2T_0 :
\|u^\lambda\|_{L^{\frac{4(n-1)}{n}}_tL^{\frac{2(n-1)}{n-2}}_x([0,t]\times
\R^n )} \leq Kt^{\frac{n-2}{4(n-1)}}\},
\end{equation}
with $K$ a constant to be chosen later. We claim that $S$ is the
whole interval $[0,\lambda^2T_0].$ Indeed, assume by contradiction
that it is not so, then since
$$\|u^\lambda\|_{L^{\frac{4(n-1)}{n}}_tL^{\frac{2(n-1)}{n-2}}_x([0,t]\times
\R^n )}$$ is a continuous function of time, there exists a time $T
\in [0,\lambda^2T_0]$ such that
\begin{align}
\label{contr1}
&\|u^\lambda\|_{L^{\frac{4(n-1)}{n}}_tL^{\frac{2(n-1)}{n-2}}_x([0,T]\times
\R^n )} > K T^{\frac{n-2}{4(n-1)}} \\
\label{contr2}
&\|u^\lambda\|_{L^{\frac{4(n-1)}{n}}_tL^{\frac{2(n-1)}{n-2}}_x([0,T]\times
\R^n )} \leq 2K T^{\frac{n-2}{4(n-1)}}.
\end{align}
From the interaction Morawetz estimate \eqref{Morfinal}, we have that
\begin{equation}\label{Morfinal2}
\|u^\lambda\|_{L^{\frac{4(n-1)}{n}}_tL^{\frac{2(n-1)}{n-2}}_x([0,T]\times
\R^n )} \lesssim T^{\frac{n-2}{4(n-1)}}
\|u_0\|_{L^2_x}^{\frac{1}{2}}\|u^\lambda\|_{L^\infty_t\dot{H}^{\frac{1}{2}}_x([0,T]\times
\R^n )}^{\frac{n-2}{n-1}}.
\end{equation}
We proceed to estimate the right-hand side of the
inequality above in terms of norms involving $Iu$ instead of $u$,
\begin{align}
\|u^\lambda(t)\|_{\dot{H}_x^{\frac{1}{2}}}
&\leq \|P_{<N}u^\lambda(t)\|_{\dot{H}_x^{\frac{1}{2}}}
+ \|P_{\geq N}u^\lambda(t)\|_{\dot{H}_x^{\frac{1}{2}}} \nonumber \\
& \leq \|P_{<N} u^\lambda(t)\|^{\frac{1}{2}}_{L_x^2}\|P_{<N} u^\lambda(t)\|^{\frac{1}{2}}_{H_x^1}
+ \frac{1}{N^{\frac{1}{2}}}\|Iu^\lambda(t)\|_{H_x^1} \label{lhint}\\
& \leq \| u_0\|^{\frac{1}{2}}_{L_x^2}\| Iu^\lambda(t)\|^{\frac{1}{2}}_{H_x^1}
+ \|I u^\lambda(t)\|_{H_x^1}, \nonumber
\end{align}
where to obtain \eqref{lhint} we used an interpolation and \eqref{mult1}.
Hence,
\begin{align}\label{Morfinal3}
&\|u^\lambda\|_{L^{\frac{4(n-1)}{n}}_tL^{\frac{2(n-1)}{n-2}}_x([0,T]\times
\R^n )} \\ &\lesssim T^{\frac{n-2}{4(n-1)}}
\|u_0\|_{L^2_x}^{\frac{1}{2}}\sup_{[0,T]}(\|
u_0\|^{\frac{1}{2}}_{L^2_x}\| I u^\lambda(t)\|^{\frac{1}{2}}_{H^1_x} +
\|I u^\lambda(t)\|_{H_x^1})^{\frac{n-2}{n-1}}\nonumber.
\end{align} We now split the interval $[0,T]$ into subintervals
$J_k$, $k=1,...,L$ in such a way
that\begin{equation}\|u^\lambda\|^{\frac{4(n-1)}{n}}_{L^{\frac{4(n-1)}{n}}_tL^{\frac{2(n-1)}{n-2}}_x(J_k\times
\R^n )}\leq \mu, \end{equation} with $\mu$ as in Proposition
\ref{lwp}. This is possible because of \eqref{contr2}. Then, the
number $L$ of possible subinterval must satisfy
\begin{equation} \label{L}
L \sim \frac{(2K)^{\frac{4(n-1)}{n}} T^{\frac{n-2}{n}}}{\mu}.
\end{equation}
From Proposition \ref{lwp} and Proposition \ref{endecay} we know that,
for any $\nu$ such that $0 < \nu \leq \min \{1,\frac{4}{n}\}s+s-1$
\begin{equation}\label{energybound}
\sup_{[0,T]}E(Iu^\lambda(t)) \lesssim E(Iu_0^\lambda) +
\frac{L}{N^{1-s +\nu}}
\end{equation}
and by our choice of $\lambda$, $E(Iu_0^\lambda)\lesssim 1.$
Hence, in order to guarantee that
\begin{equation}\label{energyb}
E(Iu^\lambda)\lesssim 1\end{equation} holds for all $t \in [0,T]$
we need to require  that
$$ L \lesssim N^{1-s+\nu}. $$
Since $T \leq \lambda^2T_0,$ according to \eqref{L}, this is fulfilled
as long as
\begin{equation} \label{LN}
\frac{(2K)^{\frac{4(n-1)}{n}} (\lambda^2
T_0)^{\frac{n-2}{n}}}{\mu} \sim N^{1-s+\nu},
\end{equation}

From our choice of $\lambda$, the expression \eqref{LN} implies that
\begin{equation}\label{Tn}
N^{1-s+\nu} \sim (2K)^{\frac{4(n-1)}{n}}
N^{\frac{2(1-s)}{s}(\frac{n-2}{n})}T_0^{\frac{n-2}{n}},
\end{equation}
where
$$0< \nu \leq \min \{1,\frac{4}{n}\}s+s-1.$$
We pick
\begin{equation} \label{nu}
\nu = \left(\frac{2(1-s)}{s}(\frac{n-2}{n})-1+s\right)+ .
\end{equation}
The choice of $\nu$ given by \eqref{nu} is permissible for $n=3$
as long as
$$\frac{2(1-s)}{s}(\frac{1}{3})-1+s <2s -1,$$
which is possible as long as $s>\frac{\sqrt{7}-1}{3}$. \\
\\

For dimension $n \geq 4$, \eqref{nu} is permissible as long as
$$\frac{2(1-s)}{s}(\frac{n-2}{n})-1+s < \frac{n+4}{n}s -1,$$
i.e.
$$s > \frac{-(n-2)+\sqrt{(n-2)^2+8(n-2)}}{4}.$$
With this choice of $\nu,$
we have that $N$ is a large number, for $T_0$ large. Then, from
\eqref{Morfinal3} and \eqref{energyb} we obtain
\begin{equation}\|u^\lambda\|_{L^{\frac{4(n-1)}{n}}_tL^{\frac{2(n-1)}{n-2}}_x([0,T]\times
\R^n )} \lesssim C T^{\frac{n-2}{4(n-1)}}\end{equation} for some
constant $C>0.$ This contradicts \eqref{contr1} for an appropriate
choice of $K$. Hence $S = [0,\lambda^2T_0]$, and $T_0$ can
be chosen arbitrarily large. In addition, we have also proved that
$$\| I u^\lambda(\lambda^{2}T_{0})\|_{H^1_x}=O(1).$$
But then,
$$\|u(T_{0})\|_{H^{s}} \lesssim \|u(T_{0})\|_{L^{2}}+\|u(T_{0})\|_{\dot{H}^{s}}=
\|u_{0}\|_{L^{2}}+\lambda^{s}\|u^{\lambda}(\lambda^{2}T_{0})\|_{\dot{H}^{s}}$$ 
\\
$$\lesssim \lambda^{s}\| I u^\lambda(\lambda^{2}T_{0})\|_{H^1_x}\lesssim \lambda^{s}\lesssim N^{1-s}
\lesssim T_{0}^{\alpha(s,n)}$$
\\
where $\alpha(s,n)$ is a positive number that depends on $s$ and $n$. Since $T_{0}$ is arbitrarily
 large, the apriori bound on the $H^{s}$ norm concludes the global well-posedness
in the range of $s$ that we summarize below.
\begin{align*}1>s_{3}&> \frac{\sqrt{7}-1}{3},\\
1>s_{n}&>\frac{-(n-2)+\sqrt{(n-2)^2+8(n-2)}}{4},\ \ n\geq 4,\\
\end{align*}
where $s_n$ denotes the Sobolev index $s$ corresponding to the space $H^s_x(\R^n)$.

%Notice that all these conditions are within the acceptable range
%$1>s>\frac{1}{1+\min\{1,\frac{4}{n}\}}$. Thus for this range of $s$
%we have that
%$$\|u^{\lambda}\|_{\dot{H}^{s}} \lesssim \|Iu^{\lambda}\|_{\dot{H}^{1}} \lesssim 1$$
%for all $t \in [0,\lambda^2T_0]$.
%This implies that
%$$\|u\|_{\dot{H}^{s}}=\lambda^{s}\|u^{\lambda}\|_{\dot{H}^{s}}
%\lesssim \lambda^{s} \lesssim N^{1-s}$$
%for any $t \in [0,T_0]$.
%Going back and expressing $N$ in terms of
%$T_{0}$ by \eqref{Tn}, we get that
%$$\alpha(s)=\frac{(n-2)s(1-s)}{(1-s)ns+\nu ns-2(n-2)(1-s)}.$$ By a
%classical density argument this shows global well-posedness and the
%polynomial bounds and concludes the proof.
\end{proof}

\end{document}